\documentclass[11pt]{amsart}
\usepackage[normalem]{ulem}
\usepackage[usenames]{color}
\usepackage{graphicx,amscd,psfrag,amssymb}
\usepackage[colorlinks=true,citecolor=red,linkcolor=blue]{hyperref}

\theoremstyle{plain}

\theoremstyle{definition}

\theoremstyle{remark}

\begin{document}

\title{A Trisectrix from a Carpenter's Square} 

\author{David Richeson}   \address{Dickinson College\\ Carlisle, PA 17013} \email{richesod@dickinson.edu} 

\begin{abstract}
In 1928 Henry Scudder described how to use a carpenter's square to trisect an angle. We use the ideas behind Scudder's technique to define a trisectrix---a curve that can be used to trisect an angle. We also describe a compass that could be used to draw the curve.
\end{abstract}

\maketitle
In 1837 Pierre Wantzel famously proved that it is impossible to trisect an arbitrary angle using only a compass and straightedge \cite{Wantzel:1837}. However, it \emph{is} possible to trisect an angle if we are allowed to add additional items to our toolkit. 

We can trisect an angle if we have a marked straightedge \cite[pp.~185--87]{Knorr:1993}, a Mira (a vertical mirror used to teach transformational geometry) \cite{Emert:1994}, a tomahawk-shaped drawing tool \cite{Sackman:1956}, origami paper \cite{Fusimi:1980}, or a clock \cite{Moser:1947}. We can also trisect an angle if we are able to use curves other than straight lines and circles: a hyperbola \cite[pp.~22--23]{Yates:1942a}, a parabola \cite[pp.~206--08]{Descartes:1954}, a quadratrix \cite[pp.~81--86]{Knorr:1993}, an Archimedean spiral \cite[p.~126]{Boyer:1991}, a conchoid \cite[pp.~20--22]{Yates:1942a}, a trisectrix of Maclaurin \cite{Reyerson:1977}, a lim\c{c}on \cite[pp.~23--25]{Yates:1942a}, and so on; such a curve is called a \emph{trisectrix}. In many cases, we can use specially-design compasses to draw these or other trisectrices. For instance, Descartes designed such a compass \cite[pp.~237--39]{Bos:2001}. 

The literature on different construction tools and techniques, new compasses, and their relationships to angle trisection and the other problems of antiquity is vast. A reader interested in learning more may begin with \cite{Bos:2001, Heath:1921a,Knorr:1993,Martin:1998,Yates:1942a}.

In this note we describe a trisection technique discovered by Henry Scudder in 1928 that uses a carpenter's square \cite{Scudder:1928}. Then we use the ideas behind this construction to produce a new trisectrix, and we describe a compass that can draw the curve.

\section{Angle Trisection Using a Carpenter's Square}
A carpenter's square---a common drawing tool found at every home improvement store---consists of two straightedges joined in a right angle. To carry out Scudder's construction we need a mark on one leg such that the distance from the corner is twice the width of the other leg. For instance, we will assume that the longer leg is one inch wide and that there is a mark two inches from the corner on the shorter leg. 

\begin{figure}[ht]
\includegraphics{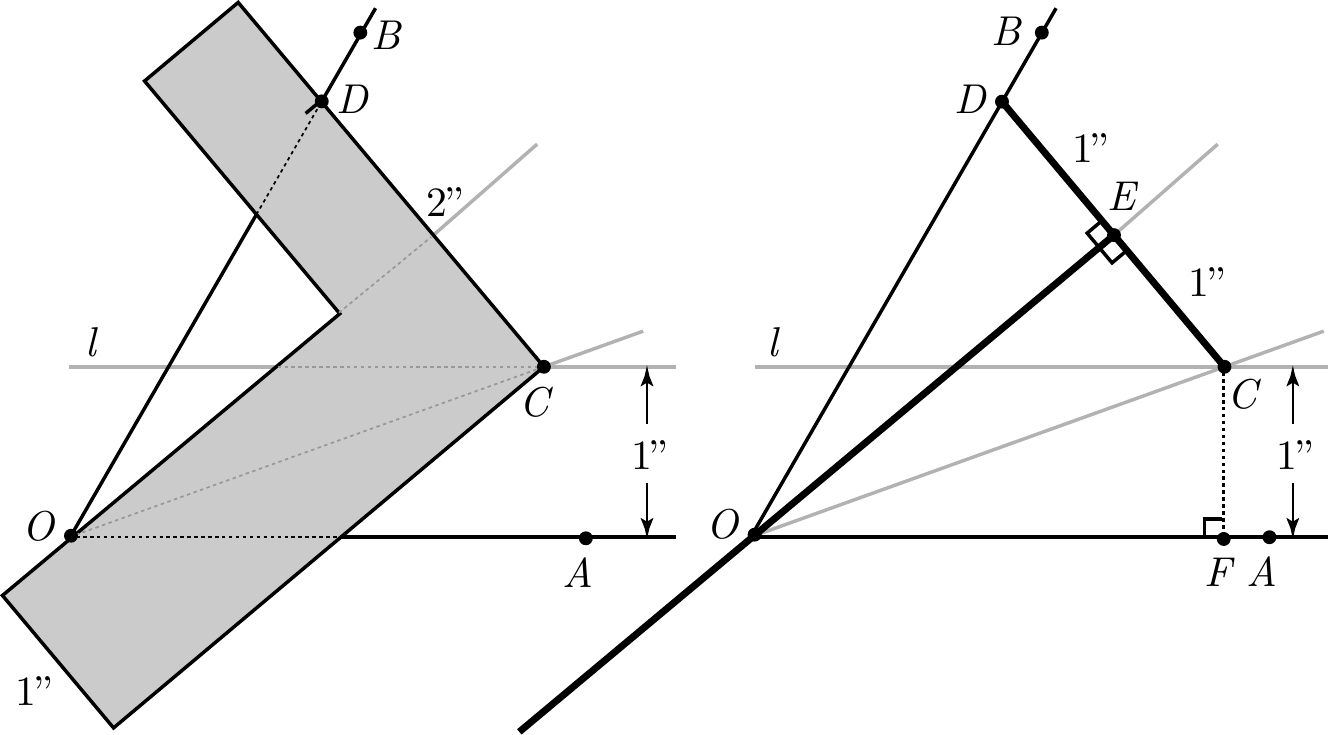}
\caption{A carpenter's square or a T-shaped tool can be used to trisect an angle.}
\label{fig:carpenterssquare}
\end{figure}

Let's say we wish to trisect the angle $\angle AOB$ in figure \ref{fig:carpenterssquare}. First, we draw a line $l$ parallel to and one inch away from $AO$; this can be accomplished using a compass and straightedge, but a simpler method is to use the long leg of the carpenter's square as a double-edged straightedge. We now perform the step that is impossible using Euclidean tools: Place the carpenter's square so that the inside edge passes through $O$, the two inch mark lies on the line $BO$, and the corner sits on the line $l$ (at the point $C$, say). Then the inside edge of the carpenter's square and the line $CO$ trisect the angle. This procedure works for any angle up to $270^{\circ}$, although the larger the angle, the narrower the short leg of the carpenter's square must be.

In fact, we do not need a carpenter's square to carry out this construction. All we need is a T-shaped device (shown on the right in figure \ref{fig:carpenterssquare}) in which the top of the T is two inches long. It is not difficult to see that this technique trisects the angle: the right triangles $COF$, $COE$, and $DOE$ in figure \ref{fig:carpenterssquare} are congruent.

\section{A New Compass}

\begin{figure}[ht]
\includegraphics{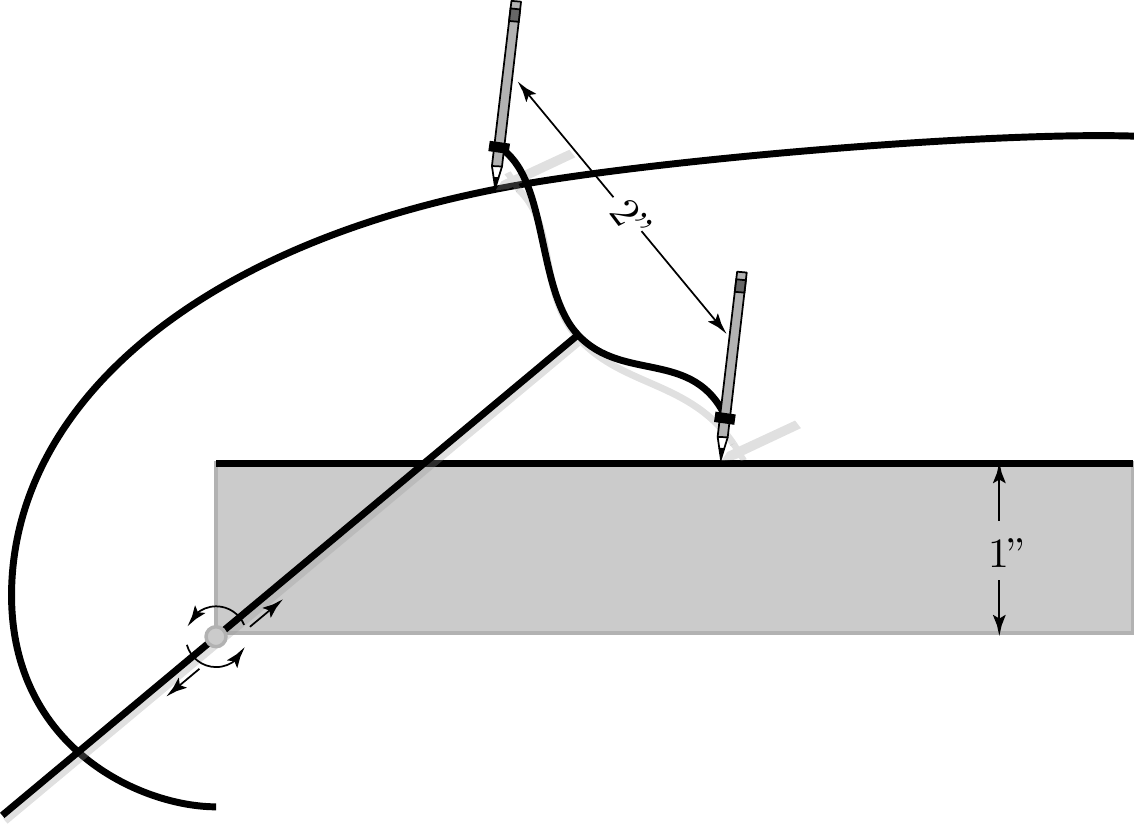}
\caption{A compass to draw the carpenter's square curve.}
\label{fig:compass}
\end{figure}

We now use the carpenter's square as inspiration to create a compass to draw a trisectrix (see figure \ref{fig:compass}). The device has a straightedge that is one inch wide and a T-shaped tool with pencils at both ends of the two-inch top of the T. The long leg of the T passes through a ring at one corner of the straightedge. The T can slide back and forth in the ring, and the ring can rotate. One pencil draws a line along the straightedge. The other pencil draws the curve we call the \emph{carpenter's square curve}.

\begin{figure}[ht]
\includegraphics{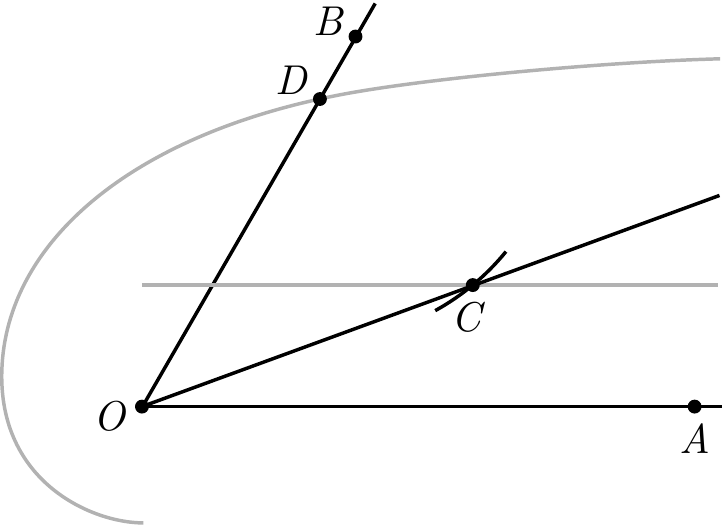}
\caption{We can use the carpenter's square curve to trisect an angle.}
\label{fig:usecurve}
\end{figure}

We use the compass as follows. Suppose we would like to trisect $\angle AOB$ in figure \ref{fig:usecurve}. Place the bottom of the straightedge along $OA$ with the ring located at $O$. Use the compass to draw the straight line $l$ and the carpenter's square curve. Say that $BO$ intersects the curve at $D$. Use an ordinary compass to draw a circle with center $D$ and a two-inch radius. It will intersect $l$ at two points. Label the right-most point (viewed from the perspective of figure \ref{fig:usecurve}) $C$. Then $OC$ trisects the angle. Use an ordinary compass and straightedge to bisect $\angle COD$ to obtain the other trisecting ray.

\section{The Carpenter's Square Curve}
What is this carpenter's square curve? Does it have a closed form? Is it algebraic or transcendental? (In \cite{Yates:1940a}, Yates used a carpenter's square in a different way to generate a different curve---a cissiod. Yates gives an algebraic expression for his curve and shows that it can be used to compute cube roots.)

\begin{figure}[ht]
\includegraphics{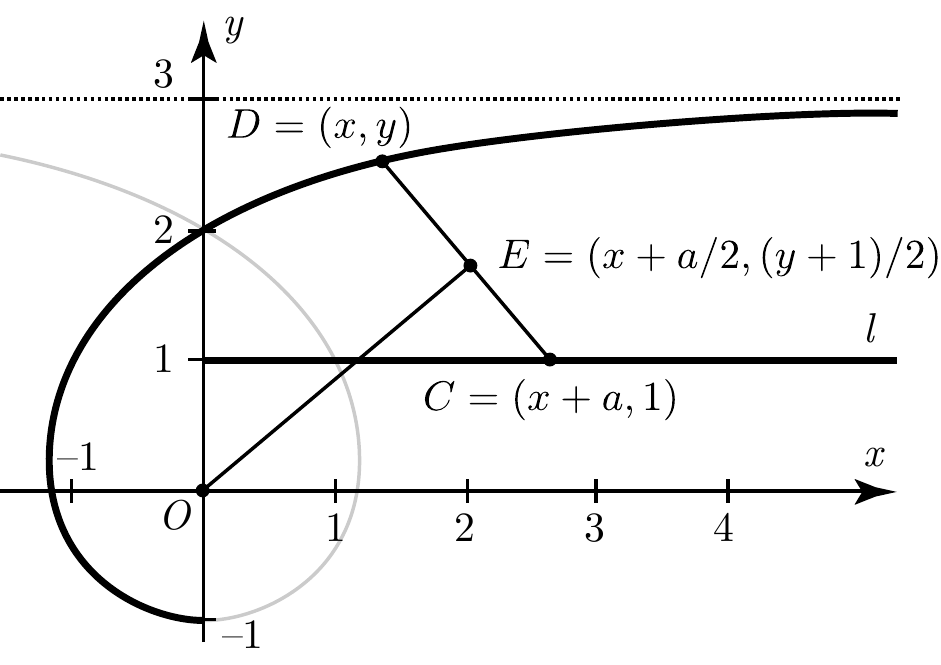}
\caption{The carpenter's square curve.}
\label{fig:curve}
\end{figure}

First we introduce $x$- and $y$-axes. Let $O$ be the origin and $OA$ be the positive $x$-axis (see figure \ref{fig:curve}). Let $D=(x,y)$ and $C=(x+a, 1)$. Because $|CD|=2$, $a^{2}+(1-y)^{2}=4$, and hence $a=\sqrt{(3-y)(y+1)}.$ (Notice that $a\ge 0$ throughout the construction.) Also, $E$ is the midpoint of $CD$, so $E=(x+a/2,(y+1)/2)$. Because $CD$ and $EO$ are perpendicular, \[\frac{(y+1)/2}{x+a/2}=-\frac{a}{1-y}.\] Substituting our expression for $a$ and simplifying, we obtain \[x^{2}=\frac{(y-2)^{2}(y+1)}{3-y}.\] This algebraic curve has a self-intersection at $(0,2),$ and $y=3$ is a horizontal asymptote. However, as we see in figure \ref{fig:curve}, our compass does not trace this entire curve.  

To see an interactive applet of this trisection, visit \emph{\href{http://ggbtu.be/mJpaNPATB}{ggbtu.be/mJpaNPATB}}.

\bibliographystyle{plain} 
\bibliography{carpenterssquare}

\begin{thebibliography}{10}

\bibitem{Bos:2001}
Henk J.~M. Bos.
\newblock {\em Redefining geometrical exactness: {D}escartes' transformation of
  the early modern concept of construction}.
\newblock Springer-Verlag, New York, 2001.

\bibitem{Boyer:1991}
Carl~B. Boyer and Uta~C. Merzbach.
\newblock {\em A History of Mathematics}.
\newblock John Wiley \& Sons, New York, 2 edition, 1991.

\bibitem{Descartes:1954}
Ren\'e Descartes.
\newblock {\em The Geometry of {R}en\'e {D}escarte: {T}ranslated from {F}rench
  and {L}atin by {D}avid {E}ugene {S}mith and {M}arcia {L}. {L}atham}.
\newblock Dover Publications Inc., New York, 1954.

\bibitem{Emert:1994}
John~W. Emert, Kay~I. Meeks, and Roger~B. Nelson.
\newblock Reflections on a {M}ira.
\newblock {\em The American Mathematical Monthly}, 101(6):544--549, 1994.

\bibitem{Fusimi:1980}
Koji Fusimi.
\newblock Trisection of angle by {A}be.
\newblock {\em Saiensu (supplement)}, page~8, October 1980.

\bibitem{Heath:1921a}
Thomas~L. Heath.
\newblock {\em A history of {G}reek mathematics. {V}ol. {I}: From {T}hales to
  {E}uclid}.
\newblock Clarendon Press, Oxford, 1921.

\bibitem{Knorr:1993}
Wilbur~Richard Knorr.
\newblock {\em The ancient tradition of geometric problems}.
\newblock Dover Publications Inc., New York, 1993.

\bibitem{Martin:1998}
George~E. Martin.
\newblock {\em Geometric constructions}.
\newblock Springer-Verlag, New York, 1998.

\bibitem{Moser:1947}
Leo Moser.
\newblock The watch as angle trisector.
\newblock {\em Scripta Math.}, 13:57, 1947.

\bibitem{Reyerson:1977}
Hardy~C. Reyerson.
\newblock Anyone can trisect an angle.
\newblock {\em Mathematics Teacher}, 70:319--321, April 1977.

\bibitem{Sackman:1956}
Bertram~S. Sackman.
\newblock The tomahawk.
\newblock {\em Mathematics Teacher}, 49:280--281, April 1956.

\bibitem{Scudder:1928}
Henry~T. Scudder.
\newblock Discussions: {H}ow to trisect an angle with a carpenter's square.
\newblock {\em The American Mathematical Monthly}, 35(5):250--251, 1928.

\bibitem{Wantzel:1837}
P.~L. Wantzel.
\newblock Recherches sur les moyens de reconna\^itre si un {P}robl\`eme de
  {G}\'eom\'etrie peut se r\'esoudre avec la r\`egle et le compas.
\newblock {\em Journal de Math\'ematiques Pures et Appliqu\'ees},
  2(1):366---372, 1837.

\bibitem{Yates:1940a}
Robert~C. Yates.
\newblock The angle ruler, the marked ruler and the carpenter's square.
\newblock {\em National Mathematics Magazine}, 15(2):61--73, 1940.

\bibitem{Yates:1942a}
Robert~C. Yates.
\newblock {\em The Trisection Problem}.
\newblock The Franklin Press, Baton Rouge, LA, 1942.

\end{thebibliography}
\end{document}